\documentclass{article}
\usepackage{amsmath}
\usepackage{parskip}
\usepackage[margin=1in]{geometry}
\usepackage{amssymb}
\usepackage{graphicx} 
\usepackage{xcolor}
\usepackage{float}
\usepackage{bm}
\usepackage{amsthm}
\usepackage{fancyhdr}
\title{On the chasing rabbit problem}

\begin{document}
\title{\textsc{Analysis of SIR Reaction diffusion system with constant birth and death rate}}
\author{\textsc{Yiting Yao}}

\maketitle
\section*{\abstractname}
\textit{This is a truncation of the second year group project at Imperial college london. In this paper, we consider a semilinear reaction diffusion system of SIR model which involves the birth rate and the death rate. We first prove the non-negativity and global existence theorem to ensure that the model makes sense. We prove the uniform convergence of the infection-free solution and study an example that separable solutions can be computed. We also focus on the steady state solution, which we prove the non-uniqueness of the solution and investigate the regularity of the general solution. In the end we also introduce an interesting phenomenon, which is called the Turing instability caused by the diffusion in the model.}

\tableofcontents
    
    \section{\textit{\textcolor{blue}{Introduction to the PDE model}}}
    
    Consider the Cauchy-Neumann problem of the SIR reaction-diffusion reaction system which is the simpler model as mentioned in [1]. Let $\Omega\subset\mathbb{R}^n$ be an open, bounded, connected set with Lipschitz boundary for some $n\in\mathbb{N}$ and $\mathcal{D}=\Omega\times\mathbb{R}^+$ be the domain, the system of PDEs is formulated as follows
\begin{align}
    \partial_tS&=\chi_S \bm{\bm{\nabla}}^2S+b-\beta SI-\nu S,\nonumber\\
    \partial_tI&=\chi_I\bm{\nabla}^2 I+\beta SI-(\gamma+\nu)I,\\
    \partial_tR&=\chi_R\bm{\nabla}^2R+\gamma I-\nu R.\nonumber
\end{align}Where $b,\gamma,\nu$ and $\beta$ are the birth rate, recovery rate, death rate and the infectious rate respectively and $\chi_S,\chi_I,\chi_R$ are the diffusion constants. The PDEs are subjected to "no flux" boundary conditions and \textit{non-negative} initial conditions \begin{align*}
   & \partial_n S=\partial_n I=\partial_n R=0\quad \text{on $(x,t)\in\partial\mathcal{D}$},\\
  & S(\bm{x},0)=S_0,\ I(\bm{x},0)=I_0,\ R(\bm{x},0)=R_0\ \text{ with }\ S_0,I_0,R_0\in C(\overline{\Omega}).  
\end{align*}Where $n$ represents the outward pointing normal on the boundary. In this work, we will assume that all the parameters are \textit{positive} unless otherwise stated. This PDE model is relatively simple, there are many extensions of the model considering other factors. For instance, in [1] and [8], fraction of compliant population and the effect of vaccination were considered.

\section{\textit{\textcolor{blue}{Global existence of the solution and the infection-free solution}}}
In this section, we are going to prove the non-negativity of the solution and then use it to show the global existence of the solution. We starts with an interesting observation: denote the total population over the domain as\[
 N(t)=\int_{\Omega}S(\bm{x},t)+I(\bm{x},t)+R(\bm{x},t)\ d\bm{x},\]Differentiating both sides with respect to $t$, we see that\[
\frac{dN}{dt}=b|\Omega|-\nu N(t)\ \Rightarrow N=N_0e^{-\nu t}+\frac{b|\Omega|}{\nu}(1-e^{-\nu t}).
\]Where we have used the divergence theorem and $N_0$ is the initial value $N(0)$, $|\Omega|$ is the Lebesgue measure of $\Omega$. This says that the growth of total population is controlled by the ratio of the birth rate to death rate. This doesn't say anything with the global existence as we will need the solution to be non-negative and no singularity. To tackle with these problems, we need a few lemmas. \vspace{4pt}

\noindent\textbf{Lemma 1 (Comparison principle):} \textit{Let $u_1,u_2\in C\left(\left[0,T\right]; C^2(\Omega)\cap C(\overline{\Omega})\right)$ that satisfy}\[
\partial_t u_1-d\bm{\nabla}^2 u_1=f_1,\quad \partial_t u_2-d\bm{\nabla}^2 u_2=f_2.
\]\textit{Where $d,T>0$, and $f$ is bounded in $\mathcal{D}$, then if $u_1\geq u_2$ on $\partial\mathcal{D}$ and $f_1\geq f_2$ in $\mathcal{D}$, then $u_1\geq u_2$ in $\mathcal{D}$.}\par

\noindent\textbf{Corollary 1:} \textit{Suppose that $u\in C\left(\left[0,T\right]; C^2(\Omega)\cap C(\overline{\Omega})\right)$ satisfies}\begin{align*}
    \partial_t u-d\bm{\nabla}^2u&=f(\bm{x},t,u)\geq0,\quad (\bm{x},t)\in\mathcal{D},\\
    u(\bm{x},0)&=u_0(\bm{x})\geq0,\quad
    \bm{x}\in\Omega,\\
    \partial_n u&=0\quad \text{on }\partial\mathcal{D},
\end{align*}\textit{then $u(\bm{x},t)\geq0$ for all $(\bm{x},t)\in\mathcal{D}$.}\par
\noindent\textsc{Proof:} This is a direct application of  lemma 1 by taking $f_2=0$ and $u_2=0$.\vspace{4pt}

\noindent\textbf{Lemma 2:} \textit{Suppose that $u\in C\left(\left[0,T\right]; C^2(\Omega)\cap C(\overline{\Omega})\right)$  such that}\begin{align*}
\partial_t u-d\bm{\nabla}^2 u&\leq f(\bm{x},t)\quad (\bm{x},t)\in\Omega\times\left[0,T\right],\\
u(\bm{x},0)&=u_0(\bm{x})\geq 0\quad x\in\Omega,\\
\partial_n u&=0\quad\forall (\bm{x},t)\in\partial\Omega\times\left[0,T\right].
\end{align*}then $\forall t\geq0$,\[
||u||_{L^\infty(\overline{\Omega})}\leq ||u_0||_{L^\infty(\overline{\Omega})}+\int_0^t||f(\bm{x},s)||_{L^\infty(\overline{\Omega})} ds.
\]\textsc{Proof:} Let $v(\bm{x},t)=||u_0||_{L^\infty(\overline{\Omega})}+\int_0^t||f(\bm{x},s)||_{L^\infty(\overline{\Omega})} ds-u(\bm{x},t)$, then we see that\begin{align*}
   & \partial_t v-d\bm{\nabla}^2 v=||f(\bm{x},t)||_{L^{\infty}(\overline{\Omega})}-\partial_t u+d\bm{\nabla}^2 u\geq0,\ &\text{(by assumption)}\\
    \Rightarrow&v(\bm{x},t)\geq 0\quad \forall (\bm{x},t)\in\Omega\times\left[0,T\right],&\text{(by corollary 1)}\\
    \Rightarrow&||u_0||_{L^\infty(\overline{\Omega})}+\int_0^t||f(\bm{x},s)||_{L^\infty(\overline{\Omega})}\ ds\geq u(\bm{x},t),\quad\forall (\bm{x},t)\in\Omega\times\left[0,T\right].
\end{align*}Hence the result holds as $u$ is uniformly continuous.\vspace{3pt}

\noindent\textbf{Proposition 1 (Non-negativity of the solution given non-negative initial data):} \textit{ Given non-negative initial data, the classical solution $(S,I,R)$ is non-negative.}\vspace{3pt}

\noindent\textsc{Proof:} We will apply corollary 1 to prove this claim. We consider each single equation. For the first equation in (1), we use the change of variable \[
\widehat{S}=e^{\beta\int_0^t I(\bm{x},s)\ ds+\nu t}S,
\]so that $\widehat{S}(\bm{x},0)=S_0\geq 0$ and on the boundary, we have\[
\partial_n \widehat{S}=e^{\beta \int_0^t I(\bm{x},s)\ ds+\nu t}\left(\beta S\int_0^t \underbrace{\partial_n I(\bm{x},s)}_{=0}\ ds+\underbrace{\partial_n S}_{=0}\right)=0.
\]And we see that\[
\widehat{S}_t-\chi_S\bm{\nabla}^2 \widehat{S}=\left(S_t+\beta IS+\nu S\right)e^{\beta \int_0^t I(\bm{x},s)\ ds+\nu t}=be^{\beta \int_0^t I(\bm{x},s)\ ds+\nu t}\geq0.
\]Hence by corollary 1, $S$ is non-negative. Similarly, to apply corollary 1 to prove the non-negativity of $I$, we define\[
\widehat{I}=e^{-\beta\int_0^t S(\bm{x},s)\ ds+(\nu+\gamma) t}I,
\]then the non-negativity follows. Finally, by letting $\widehat{R}=e^{\nu t}R$, we can conclude the non-negativity of $R$. \vspace{3pt}

\noindent\textbf{Theorem 1 (Global existence of the solution)} \textit{Let $S$, $I$ and $R$ be the classical solution to (1) on $\left[0,T^*\right)$, where $T^*$ is defined to be}\[ T^*:=\sup\{\tau>0:\text{there exists a classical solution to $(1)$}\}
\]with non-negative initial data $S_0,I_0,R_0\in L^\infty(\Omega)$. If $T^*<\infty$, then there exists a constant $M=M(T^*)$ such that \[
||S||_{L^\infty(\overline{\Omega})},||I||_{L^\infty(\overline{\Omega})},||R||_{L^\infty(\overline{\Omega})}\leq M(T^*)
\]\textsc{Proof:} By proposition 1, $S,I,R$ are always non-negative given non-negative initial data, hence we see that\[
\partial_tS-\chi_S\bm{\nabla}^2S\leq b.
\]By lemma 2, \[
||S(\bm{x},t)||_{L^\infty(\overline{\Omega})}\leq ||S_0||_{L^\infty(\overline{\Omega})}+T^*||b||_{L^\infty(\overline{\Omega})}.
\]Hence $S(\bm{x},t)$ is bounded, using the similar method to bound $I$ and $R$, we see that the claim holds.\vspace{2pt}

\noindent\textbf{Remark:} To prove the local/global existence, an alternative way to do this is to use a method called \textit{invariant region method}, the idea is to prove $F\cdot\hat{\bm{n}}<0$ on $\partial S$ for some $S\subset\Omega$, where $\hat{\bm{n}}$ is the outward pointing normal and $F$ is the reaction terms on the RHS of the system. The power of this method can be seen in the following parts. \vspace{2pt}

Suppose that we want to find the \textit{infection-free solution}, i.e. the solution when $\beta\ll1$. If we can prove the uniform convergence of the solution with $\beta\rightarrow0$ to the solution when $\beta=0$.  Then we can find an analytic solution as the PDEs are then linear. To start with the proof of uniform convergence, we will first identify a invariant set that can be controlled by $\beta$ that helps our further analysis when 
$\beta$ is small. \vspace{2pt}

\noindent\textbf{Lemma 3 (Invariant region):} \textit{Given non-negative initial data, in the case $\beta\ll1$}, the region \[
\{(S,I,R)\in\mathbb{R}^3:S\in\left[0,c_1\right],\ I\in\left[0,c_2\right],\ R\in\left[0,c_3\right]\}
\]is an invariant region with \begin{align*}
c_1=\frac{\gamma+\nu}{\sqrt[3]{\beta}},\quad c_2=\frac{b}{\sqrt[3]{\beta}(\gamma+\nu)},\quad c_3=\frac{\gamma}{\sqrt[3]{\beta}\nu(\gamma+\nu)}+1.
\end{align*}\textsc{Proof}: Suppose that we don't know $c_1,c_2$ and $c_3$. We use the invariant region method. And there are six boundaries we need to look at. First we rewrite the system of PDE as\[
\partial_t\bm{U}=D\bm{U}+\bm{F}(\bm{U})\quad \text{in }\mathcal{D},
\]where $D=$diag$(\chi_S,\chi_I,\chi_R)$ is the diagonal matrix and $\bm{F}(\bm{U})$ is the RHS of the equation, i.e. $(b-\beta SI-\nu S,\ \beta SI-(\gamma+\nu)I,\ \gamma I-\nu R)^T$. As explained above, we need $\bm{F}(\bm{U})\cdot\hat{\bm{n}}<0$, so we need to check the six cases separately. For $S=0$, we check the boundary behaviour\[
\bm{\nabla}_{S,I,R}(-S)\cdot\bm{F}|_{S=0}=-b<0.
\]When $S=c_1$,  we have\[
\bm{\nabla}_{S,I,R}(S)\cdot\bm{F}|_{S=c_1}=b-\beta c_1 I-\nu c_1<0.
\]By picking $c_1>b/(\beta c_2+\nu)$, then we are done, and $c_1$ will be determined later. Next, when $I=0$, we see that\[
\bm{\nabla}_{S,I,R}(-I)\cdot\bm{F}|_{I=0}=0,
\]This doesn't say anything, but we see that if $I$ approaches to $0$ from negative value, this product is negative by picking $c_1<(\gamma+\nu)/\beta$. In the case $I=c_2$, we see that\[
\bm{\nabla}_{S,I,R}(I)\cdot\bm{F}|_{I=c_2}=c_2(\beta c_1-(\gamma+\nu))<0,
\]which holds by the same condition as above. Similarly, we can find the condition on $c_3$, which is $c_3>\gamma c_2/\nu$. Combining these inequalities, we can pick (in order) \[
c_2=\frac{b}{\sqrt[3]{\beta}(\gamma+\nu)},\quad c_1=\frac{\gamma+\nu}{\sqrt[3]{\beta}},\quad c_3=\frac{\gamma}{\nu}c_2+1.
\]
\noindent\textbf{Remark:} We say controlled by $\beta$ means that \[
\lim_{\beta\rightarrow0} \beta c_1=\lim_{\beta\rightarrow0} \beta c_2=\lim_{\beta\rightarrow0} \beta c_1c_2=0.
\]
\noindent Next, we can use this lemma to prove the uniform convergence of the solution.\vspace{2pt}

\noindent\textbf{Proposition 2 (Uniform convergence of the solution):} \textit{Let $(S_\beta,I_\beta, R_\beta)$ and $(S_0,I_0,R_0)$ to be the classical solutions to (1) and the following system respectively}\begin{align}
    \partial_tS&=\chi_S \bm{\nabla}^2S+b-\nu S,\nonumber\\
    \partial_tI&=\chi_I\bm{\nabla}^2 I-(\gamma+\nu)I,\\
    \partial_tR&=\chi_R\bm{\nabla}^2R+\gamma I-\nu R.\nonumber
\end{align}\textit{With the same boundary and non-negative initial data, we have $(S_\beta,I_\beta, R_\beta)$ converging uniformly to $(S_0,I_0,R_0)$ as $\beta\rightarrow0$, i.e.}\[
\limsup_{\beta\rightarrow0}|S-S_\beta|=\limsup_{\beta\rightarrow0}|I-I_\beta|=\limsup_{\beta\rightarrow0}|R-R_\beta|=0.
\]

\noindent\textsc{Proof:} The proof is similar as above, we take $A=S_\beta-S_0$, $B=I_\beta-I_0$ and $C=C_\beta-C_0$, then we see that\[
\partial_tA-\chi_S\bm{\nabla}^2 A+\nu A=\beta S_\beta I_\beta,\quad \partial_tB-\chi_I\bm{\nabla}^2 B+(\gamma+\nu) B=-\beta S_\beta I_\beta,\quad \partial_tC-\chi_R\bm{\nabla}^2 C-\gamma B+\nu C=0.
\]With zero initial data and no flux boundary condition. Now we use the change of variables $\widehat{A}=e^{\nu t}A$, $\widehat{B}=e^{(\gamma+\nu)t}B$ and $\widehat{C}=e^{\nu t}C$, then we can rewrite our equation as\[
\partial_t \widehat{A}-\chi_S \widehat{A}=\beta S_\beta I_\beta e^{\nu t},\quad \partial_t\widehat{B}-\chi_I\bm{\nabla}^2\widehat{B}=-\beta S_\beta I_\beta e^{(\gamma+\nu)t},\quad\partial_t\widehat{C}-\chi_R\bm{\nabla}^2\widehat{C}=e^{(\nu-\gamma)t}B.
\]Then by lemma 2 (if we replace $f(\bm{x},t)$ by $f(\bm{x},t,u)$, then the inequality is obviously still correct), we can bound $\widehat{A}$,\[
||\widehat{A}||_{L^\infty(\overline{\Omega})}\leq\underbrace{||\widehat{A}(0)||_{L^\infty(\overline{\Omega})}}_{=0}+\int_0^t||\beta S_\beta I_\beta e^{\nu s}||_{L^{\infty}(\overline{\Omega})}\ ds.
\]By lemma 3, the $L^\infty$ norm of $S_\beta$ and $I_\beta$ are no greater than $c_1$ and $c_2$, so we see that \[
||A||_{L^\infty(\overline{\Omega})}\leq c_1c_2\beta.
\] Letting $\beta\rightarrow0$, we see that $S_\beta$ converges to $S_0$ uniformly, the other two uniform convergence result can be obtained by using the same method.\vspace{2pt}

\noindent\textbf{\textit{Example}:} Now we can find the infectious-free solution. Notice that the resultant PDEs are now linear and it turns out that $S$ and $I$ doesn't affect each other. Now we consider the case when $\Omega=\left(0,L\right)^2$ for some $L>0$. We can simplify the PDEs by using the change of variables \[
\overline{S}=\left(-\frac{b}{\nu}+S\right)e^{\nu t},\quad \overline{I}=e^{(\gamma+\nu)t}I,\quad \overline{R}=e^{\nu t}R.
\]Then the equations now become\[
\partial_t\overline{S}=\chi_S\bm{\nabla}^2\overline{S},\quad \partial_t\overline{I}=\chi_I\bm{\nabla}^2\overline{I},\quad \partial_t\overline{R}=\chi_R\bm{\nabla}^2\overline{R}+\gamma e^{-\gamma t}\overline{I}.
\]The solution can be obtained using separation of variables. With some efforts, we see that\begin{align*}
\overline{S}=\sum_{n=0}^{\infty}\sum_{m=0}^\infty &c_{n,m}\cos\left(\frac{n\pi}{L}x\right)\cos\left(\frac{n\pi}{L}y\right)e^{-\frac{\chi_S\pi^2}{L^2}(m^2+n^2)t},\\
\overline{I}=\sum_{n=0}^{\infty}\sum_{m=0}^\infty &d_{n,m}\cos\left(\frac{n\pi}{L}x\right)\cos\left(\frac{n\pi}{L}y\right)e^{-\frac{\chi_I\pi^2}{L^2}(m^2+n^2)t},\\
\overline{R}=\sum_{n=0}^{\infty}\sum_{m=0}^\infty &e_{n,m}\cos\left(\frac{n\pi}{L}x\right)\cos\left(\frac{n\pi}{L}y\right)e^{-\frac{\chi_I\pi^2}{L^2}(m^2+n^2)t-\gamma t}\\+&f_{n,m}\cos\left(\frac{n\pi}{L}x\right)\cos\left(\frac{n\pi}{L}y\right)e^{-\frac{\chi_R\pi^2}{L^2}(m^2+n^2)t}. 
\end{align*}
Where $\overline{R}$ is guaranteed by uniqueness of solution. $c_{n,m},d_{n,m},e_{n,m}$ and $f_{n,m}$ can be computed using the initial data.\begin{align*}
c_{n,m}&=\frac{4}{L^2}\int_{\left[0,L\right]^2}S_0\cos\left(\frac{n\pi}{L}x\right)\cos\left(\frac{n\pi}{L}y\right)\ dxdy,\\
d_{n,m}&=\frac{4}{L^2}\int_{\left[0,L\right]^2}I_0\cos\left(\frac{n\pi}{L}x\right)\cos\left(\frac{n\pi}{L}y\right)\ dxdy,\\
e_{n,m}&=\frac{\gamma L^2 d_{n,m}}{\gamma L^2+\pi^2(m^2+n^2)(\chi_I-\chi_R)},\\
f_{n,m}&=-e_{n,m}+\frac{4}{L^2}\int_{\left[0,L\right]^2}R_0\cos\left(\frac{n\pi}{L}x\right)\cos\left(\frac{n\pi}{L}y\right)\ dxdy.
\end{align*}$S,I,R$ can be computed by substituting the variables back. We let the initial data to be a Gaussian-like distribution\[
S_0=1-e^{-(x-2.5)^2-(y-2.5)^2},\quad I_0=\frac{e^{-(x-2.5)^2-(y-2.5)^2}}{2},\quad R_0=\frac{e^{-(x-2.5)^2-(y-2.5)^2}}{2},
\]with $\chi_S=0.3$, $\chi_I=0.4$, $\chi_R=0.5$, $b=0.5$, $\beta=0.01$, $\nu=0.5$, $\gamma=0.5$ and $L=5$. The visualisation of $(S_\beta,I_\beta,R_\beta)$ and $(S,I,R)$ can be seen in figure 1.
\begin{figure}[H]
    \centering
    \begin{minipage}[b]{0.14\textwidth}
        \centering
        \includegraphics[width=\textwidth]{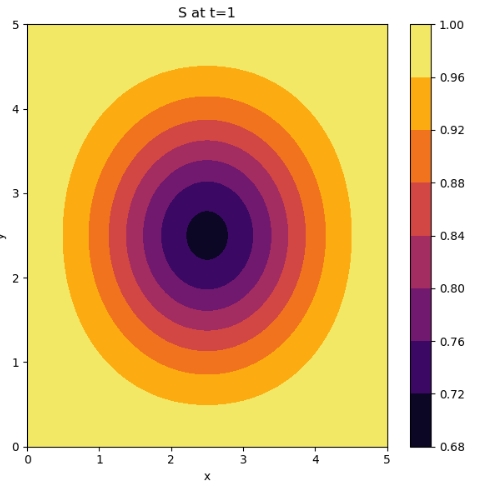}
    \end{minipage}
    \hspace{0.1cm}
    \begin{minipage}[b]{0.14\textwidth}
        \centering
        \includegraphics[width=\textwidth]{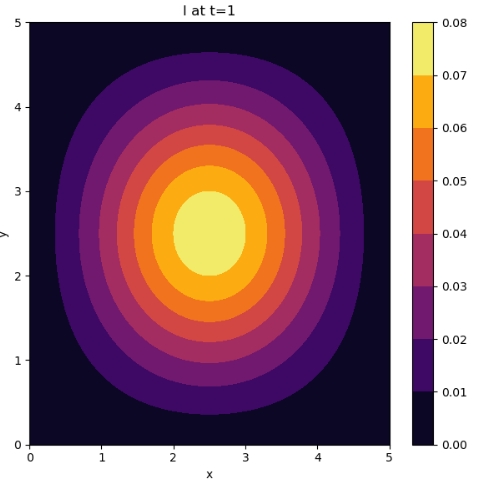}
    \end{minipage}
    \hspace{0.1cm}
    \begin{minipage}[b]{0.14\textwidth}
        \centering
        \includegraphics[width=\textwidth]{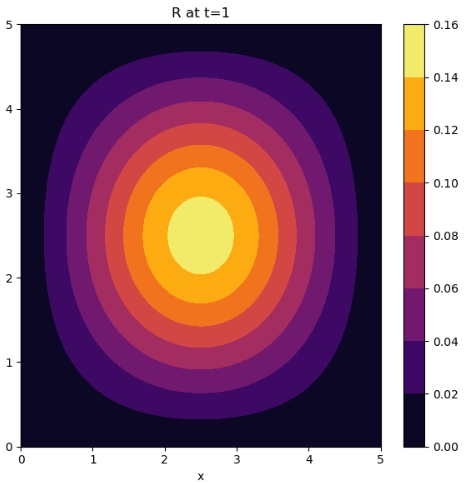}
    \end{minipage}
    \hspace{0.1cm}
    \begin{minipage}[b]{0.14\textwidth}
        \centering
        \includegraphics[width=\textwidth]{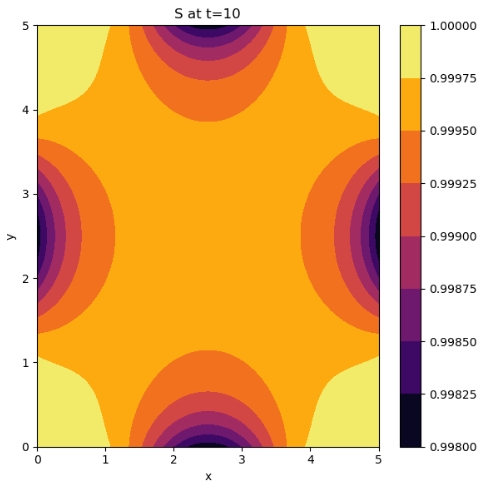}
    \end{minipage}
    \hspace{0.1cm}
    \begin{minipage}[b]{0.14\textwidth}
        \centering
        \includegraphics[width=\textwidth]{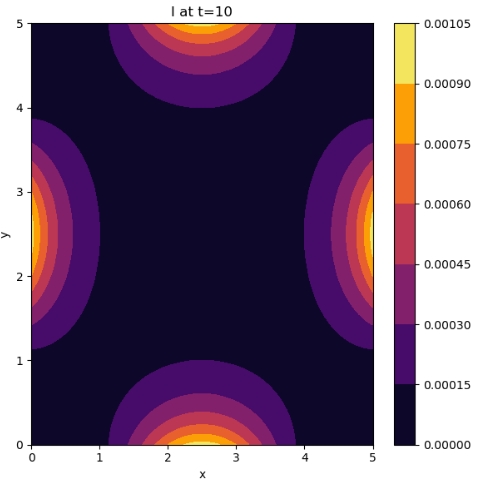}
    \end{minipage}
    \hspace{0.1cm}
    \begin{minipage}[b]{0.14\textwidth}
        \centering
        \includegraphics[width=\textwidth]{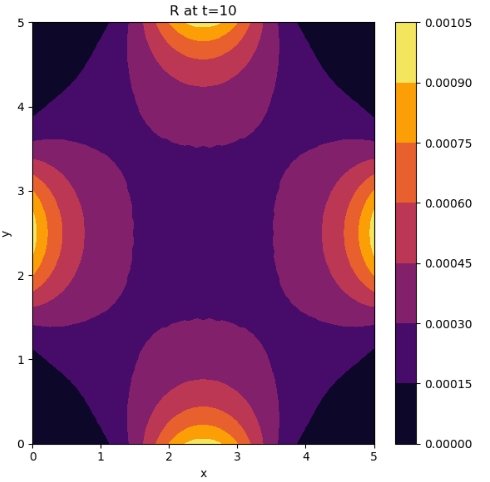}
    \end{minipage}
    \caption{Visualisation of $(S_\beta,I_\beta,R_\beta)$ at time $t=1$ and $t=10$.}
\end{figure}

\noindent\textbf{Remark:}  A sufficient but not necessary condition for the uniform convergence of the Fourier Series is that $S,I,R\in L^2(\Omega)$. Besides, as $\beta$ stands for the infectious rate, so we are expecting all the infectious die out as time goes large, and the plot of $I$ indicated this. We can also verify this analytically, by multiplying $I$ on both sides of the second equation in (2), we see that\[
\frac{1}{2}\frac{d}{dt}||I||_{L^2}^2\leq -(\gamma+\nu)||I||_{L^2}^2\Rightarrow ||I||_{L^2}^2\leq Ce^{-2(\gamma+\nu)t}\rightarrow0\quad\text{as $t\rightarrow\infty$.}
\]Using the H\"older inequality, one can also show that $||R||_{L^2}\rightarrow0$ as $t\rightarrow\infty$. As we picked the birth rate and the death rate to be the same, so $S\rightarrow1$.
\section{\textit{\textcolor{blue}{A Priori Estimates on the steady state solution}}}
Our further analysis will be on the special solution to the parabolic system-the steady state solution, which is independent of the time variable, i.e. the solution of\begin{align}
    -\chi_S \bm{\nabla}^2S&=b-\beta SI-\nu S,\nonumber\\
    -\chi_I\bm{\nabla}^2 I&=\beta SI-(\gamma+\nu)I,\\
    -\chi_R\bm{\nabla}^2R&=\gamma I-\nu R.\nonumber
\end{align}Subjecting to the Neumann boundary condition\begin{align*}
    \partial_n S=\partial_n I=\partial_n R=0\quad \text{on $(x,t)\in\partial\Omega$}.
\end{align*}
In this section, we will do some A Priori estimates on the classical solution. Again, if the solution exists, we may have the following equation: (1) Under the assumption that the solution is \textit{non-negative}, is the solution unique? (2) What is the regularity of the general solution under some assumptions? To do so, we need a comparison principle on elliptic equation, as well as a uniqueness lemma on linear PDE.\vspace{3pt}

\noindent\textbf{Theorem 2:} \textit{Let $g\in C^1(\overline{\Omega}\times\mathbb{R})$.}\begin{itemize}
    \item If $w\in C^2(\Omega)\cap C^1(\overline{\Omega})$ satisfies\[
    \bm{\nabla}^2 w+g(x,w)\geq 0\quad \text{in $\Omega$},\quad \partial_n w\leq 0\quad \text{on }\partial\Omega,
    \]and $w(x_0)=\max_{\overline{\Omega}}w$, then $g(x_0,w(x_0))\geq0$.
    \item If $w\in C^2(\Omega)\cap C^1(\overline{\Omega})$ satisfies\[
    \bm{\nabla}^2 w+g(x,w)\leq 0\quad \text{in $\Omega$},\quad \partial_n w\geq 0\quad \text{on }\partial\Omega,
    \]and $w(x_0)=\min_{\overline{\Omega}}w$, then $g(x_0,w(x_0))\leq0$. 
 \end{itemize}

\noindent\textbf{Lemma 4}: \textit{Let $a>0$ and $b\in\mathbb{R}$, then the classical solution $u\in C^2(\Omega)$ to the following linear PDE, subjecting to the Neumann boundary condition $\partial_nu=0$ on $\partial\Omega$, is unique.}\begin{equation}
\bm{\nabla}^2 u=au+b\quad \text{in $\Omega$}.
\end{equation}\textsc{Proof:} Assume for non-uniqueness of the solution, let $u_1$ and $u_2$ be two distinct solutions to (9). Define $w=u_1-u_2$, then $w$ is the solution to the equation\[
\bm{\nabla}^2w=aw\quad \text{in $\Omega$},
\]with no flux boundary condition. Hence using the divergence theorem, we see that\[
\int_{\Omega} w\bm{\nabla}^2 w+|\bm{\nabla} w|^2\ dV=\underbrace{\int_{\partial\Omega}w\partial_nw\ dS}_{=0}=\int_{\Omega} aw^2+|\bm{\nabla} w|^2\ dV,
\]since $a>0$, we  see that $w=0$ in $\Omega$, so the uniqueness follows.\vspace{3pt}

\noindent To establish the uniqueness result, we first note that the equilibria of the system (3) are found to be \begin{equation}
\boxed{\bm{A}^*_1=\left(\frac{b}{\nu},0,0\right),\quad\bm{A}^*_2=\left(\frac{\gamma+\nu}{\beta},\frac{b\beta-\nu\gamma-\nu^2}{\beta(\gamma+\nu)},\frac{\gamma(b\beta-\nu\gamma-\nu^2)}{\nu\beta(\gamma+\nu)}\right).}
\end{equation}Which are obtained by computing the roots of RHS of (4), i.e. \[
b-\beta SI-\nu S=0,\quad \beta SI-(\gamma+\nu)I=0,\quad \gamma I-\nu R=0.
\]
\noindent\textbf{Proposition 3:} \textit{Under the assumption that the classical solution to the elliptic system (3) is non-negative, then}\begin{itemize}
    \item \textit{If $b\beta-\nu\gamma-\nu^2<0$, i.e. $\bm{A}^*_2$ doesn't exist, then the system (3) has a unique classical solution $\bm{A}^*_1$.}
    \item \textit{If $b\beta-\nu\gamma-\nu^2\geq 0$, i.e. $\bm{A}^*_2$ exists, then the classical solution to (3) is either $\bm{A}^*_1$ or $\bm{A}^*_2$. }
\end{itemize}\vspace{3pt}

\noindent\textsc{Proof:} Let the classical solution to (3) to be $(S,I,R)$. Let $x_1$ and $x_2$ be the maximum points of $S$ and $I$ respectively, and $x_3$ and $x_4$ to be the minimum points of $S$ and $I$ respectively, by theorem 2, we have\begin{align}
    b-\beta S(x_1)I(x_1)-\nu S(x_1)&\geq0,\\
    b-\beta S(x_3)I(x_3)-\nu S(x_3)&\leq0,\\
    \beta S(x_2) I(x_2)-(\gamma+\nu) I(x_2)&\geq0,\\
    \beta S(x_4) I(x_4)-(\gamma+\nu) I(x_4)&\leq0.
\end{align}From (8), we see that either $(a)$: $I(x_2)$=0 or $(b)$: $S(x_2)\geq\frac{\gamma+\nu}{\beta}$. If it's the case that $(a)$ is true, then we see that $0\leq I(x)\leq 0$, so $I=0$, thus by lemma 4, it follows that $S=b/v$ and $R=0$ is the unique solution, which means that the solution is constant and equal to $A^*_1$. If $(b)$ is the case, we simply have \begin{equation}
    S(x_1)\geq S(x_2)\geq\frac{\gamma+\nu}{\beta}
\end{equation}
So, from (6), we can deduce that\begin{equation}
I(x_1)\leq\frac{b}{\beta S(x_1)}-\frac{\nu}{\beta}\leq \frac{b\beta-\gamma\nu-\nu^2}{\beta(\gamma+\nu)}\Rightarrow I(x_4)\leq\frac{b\beta-\gamma\nu-\nu^2}{\beta(\gamma+\nu)}.
\end{equation}Similarly, we can obtain the lower bound of $I(x_2)$ by using (7) and (9), which is\[
I(x_2)\geq \frac{b\beta-\gamma\nu-\nu^2}{\beta(\gamma+\nu)}.
\]Hence $I$ is solved by the constant solution $\frac{b\beta-\gamma\nu-\nu^2}{\beta(\gamma+\nu)}$, and then it follows again from lemma 4 that the solution to the system is uniquely determined to be $\bm{A}^*_2$ provided it exists, i.e. $b\beta-\nu\gamma-\nu^2\geq 0$.\vspace{3pt}

Clearly, if the non-negativity assumption is dropped, there are many other possible solutions, an example is shown in figure 2.
\begin{figure}[H]
    \centering
    \begin{minipage}[b]{0.3\textwidth}
        \centering
        \includegraphics[width=\textwidth]{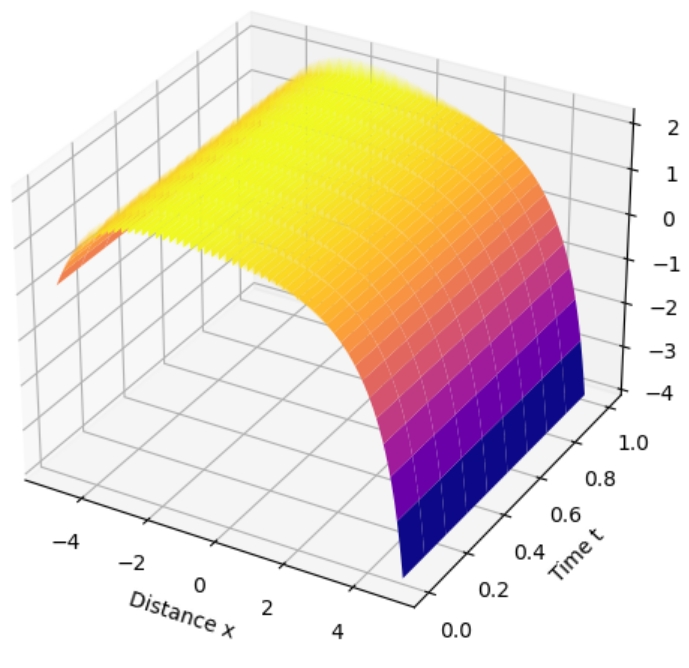}
        
    \end{minipage}
    \hspace{0.1cm}
    \begin{minipage}[b]{0.3\textwidth}
        \centering
        \includegraphics[width=\textwidth]{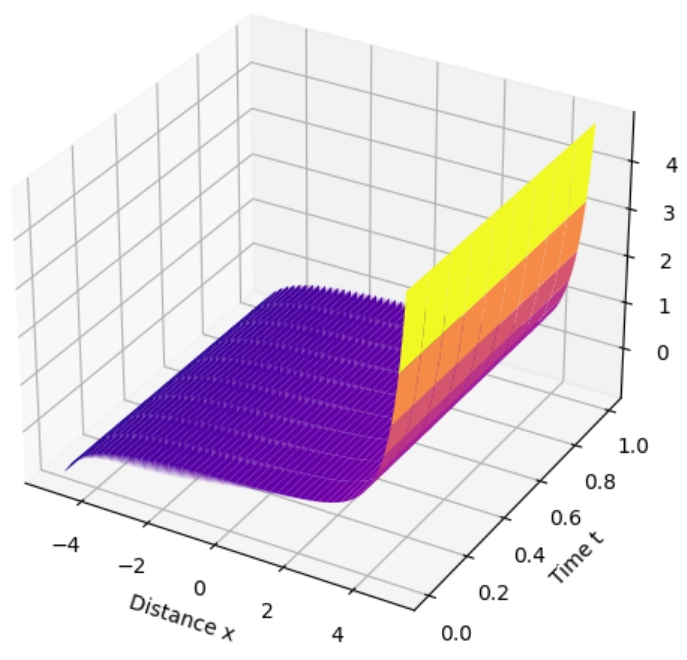}
    
    \end{minipage}
    \hspace{0.1cm}
    \begin{minipage}[b]{0.3\textwidth}
        \centering
        \includegraphics[width=\textwidth]{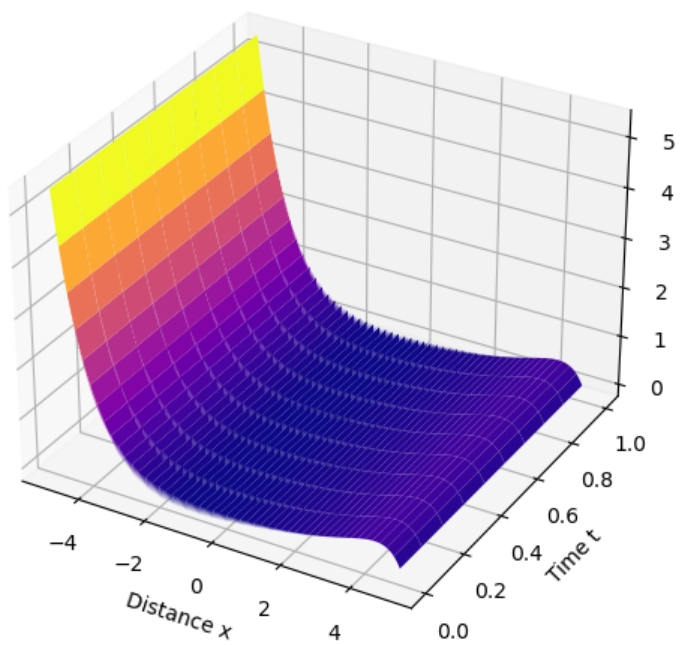}
     
    \end{minipage}
    \caption{The time independent one-dimensional steady state solution on $(-5,5)$.}
\end{figure}

 In fact, the non-negativity restricts many possibilities of the solution, proposition 3 has already indicated the non-uniqueness of the classical solution. Now we will give an estimate on the regularity of the \textit{general} solution.\vspace{3pt}

\noindent\textbf{Proposition 4:} \textit{If the classical solution $S,I,R\in C^2(\Omega)\cap L^3(\Omega)$, then $S,I,R\in H^1(\Omega)$}.\vspace{3pt}

\noindent\textsc{Proof:} Clearly $S,I,R\in L^2(\Omega)$ in this bounded domain. We first prove that $||R||_{H^1(\Omega)}$ can be bounded by $||I||_{L^2(\Omega)}$. Multiplying $R$ on both sides of the third equation of (8) and integrating over $\Omega$ yields,\[
-\chi_R\int_\Omega\bm{\nabla}^2 R=\gamma\int_{\Omega}IR-\nu \int_\Omega R^2,
\]by the divergence theorem, together with the H\"older inequality, we see that\[
\chi_R\int_\Omega|\bm{\nabla} R|^2\leq \gamma\left(\int_\Omega I^2\int_\Omega R^2\right)^{\frac{1}{2}}-\nu\int_\Omega R^2.
\]Which is regarded as a quadratic, so we can find the upper bound of $\int_\Omega R^2$ and $\int_\Omega |\bm{\nabla} R|^2$\[
\int_\Omega |\bm{\nabla} R|^2\leq \frac{\gamma^2}{4\nu}\int_{\Omega} I^2,\quad\text{and}\quad 
\gamma\left(\int_\Omega I^2\right)^{\frac{1}{2}} -\nu\left(\int_\Omega R^2\right)^\frac{1}{2}\geq 0\ \Rightarrow\int_\Omega R^2\leq \frac{\gamma^2}{\nu ^2}\int_\Omega I^2.
\]Now we are left to bound $\int_\Omega |\bm{\nabla} S|^2$ and $\int_\Omega |\bm{\nabla} I|^2$, summing the first two equation and multiplying by $S$ and then integrate over $\Omega$, we see that\begin{equation}
\chi_S\int_\Omega |\bm{\nabla} S|^2+\chi_I\int_\Omega \bm{\nabla} I\cdot\bm{\nabla} S=b\int_\Omega S-\nu\int_\Omega S^2-(\gamma+\nu)\int_\Omega SI.
\end{equation}It turns out we need to bound $\int_\Omega SI$ first. By integrating the second equation over $\Omega$, we see that\[
\int_\Omega SI=\frac{\gamma+\nu}{\beta}\int_\Omega I\leq \frac{\gamma+\nu}{\beta}|\Omega|^\frac{1}{2}\left(\int_\Omega I^3\right)^\frac{1}{2},
\]Thus we can then bound (12), with the use of AM-GM inequality separating $\bm{\nabla} S\cdot\bm{\nabla} I$ and the bound on $\int_\Omega SI$,\[
\chi_S\int_\Omega|\bm{\nabla} S|^2\leq \frac{\chi_I}{2}n^2\int_\Omega|\bm{\nabla} S|^2+\frac{\chi_I}{2n^2}\int_\Omega|\bm{\nabla} I|^2+b\left(\int_\Omega S^2\right)^\frac{1}{2}|\Omega|^\frac{1}{2}+\nu\int_\Omega S^2+\frac{\gamma+\nu}{\beta}|\Omega|^\frac{1}{2}\left(\int_\Omega I^2\right)^\frac{1}{2}.
\]Where $n\neq 0$. If $2\chi_S-\chi_I n^2>0$, we can simply write this as\[
\int_\Omega|\bm{\nabla} S|^2\leq\frac{\chi_I}{n^2(2\chi_S-\chi_I n^2)}\int_\Omega|\bm{\nabla} I|^2+ F(\gamma,\nu,\beta,|\Omega|,n,\chi_I,\chi_S,||I||_{L^2},||S||_{L^2}).
\]Where $F$ is an unimportant bounded function. Finally we estimate $||\bm{\nabla}I||_{L^2(\Omega)}$ by multiplying $I$ to the second equation in (3) and integrating directly,
\[
\chi_I\int_\Omega|\bm{\nabla}I|^2=\beta\int_\Omega SI^2-(\gamma+\nu)\int_\Omega I^2.
\]Using H\"older's inequality again, we see that\[
\int_\Omega SI^2\leq\left(\int_\Omega|S|^3\right)^{\frac{1}{3}}\left(\int_\Omega|I|^3\right)^\frac{2}{3}
\]which is bounded, therefore we see that $S,I,R\in H^1(\Omega)$.

Finally, we will prove that, under the same assumption as above, the solution of the elliptic system (3) is always smooth (infinitely differentiable). For simplicity, we will consider the cases when $n=1,2$ and $3$, which are more practical values.\vspace{5pt}

\noindent\textbf{Theorem 11} Under the same assumption as in Proposition 4, for $n=1,2,3$, $\forall k\in\mathbb{N}^+$, we have $S,I,R\in C^{k}(\Omega')$ for all $\Omega'\subset\subset\Omega$. ($\Omega'\subset\subset\Omega$ means that $\exists \mathcal{K}\subset\Omega$ compact and $\Omega'\subset\mathcal{K}\subset \Omega$.)\vspace{5pt}

\noindent\textsc{Proof:} We look at each case separately. In the first case $n=1$, recall that a Sobolev embedding result says that $W^{1,2}(\Omega)\hookrightarrow C^{0,\alpha}(\Omega)$, where $\alpha=1/2$. Under the assumptions of Proposition 2.34, $S,I,R\in W^{1,2}(\Omega)=H^1(\Omega)$, So $S,I,R\in C^{0,\alpha}(\Omega)$. Hence since the RHS of (3) are in $C^{0,\alpha}(\Omega)$, by the result (theorem 6.17 from [5]), we see that $S,I,R\in C^{2,\alpha}(\Omega)$. Hence the statement holds by repeatedly applying this theorem. (i.e. $u\in C^{k,\alpha}(\Omega)\Rightarrow u\in C^{k+2,\alpha}(\Omega)$.)\vspace{2pt}

In the case $n=2$, (we cannot use the Sobolev embedding result) we first notice that the RHS of (3) are clearly bounded, so by theorem 2.30 in [6] (which says $H^1(\Omega)\hookrightarrow L^q(\Omega)$ for any $q\geq 1$) and the Interior regularity Theorem (Theorem 4.10.1 in [7]), we see that $S,I,R\in W^{2,2}(\Omega')$ (We need $\bm{\nabla}^2S,\bm{\nabla}^2I,\bm{\nabla}^2R\in L^2(\Omega)$ to apply the theorem, so we need to bound the $L^2$ norm of the only non-linear term $\int_\Omega S^2I^2\leq\left(||S||_{L^4}||I||_{L^4}\right)^{1/2}$). Again, by the Sobolev embedding theory, we have $W^{2,2}(\Omega')\hookrightarrow C^{0,\alpha}(\Omega')$, thus the arguments holds by applying theorem 6.17 from [5].\vspace{2pt}

In the case $n=3$, the Gagliardo-Nirenberg-Sobolev Inequality implies that $H^1(\Omega)\hookrightarrow L^6(\Omega)$, hence we apply the H\"older inequality twice on $\int_\Omega S^2I^2$ twice to see that \[
\int_\Omega I^3\leq\left(\int_\Omega I^6\right)^\frac{1}{2}|\Omega|^{\frac{1}{2}}\Rightarrow  \int_\Omega S^2I^2\leq \left(\int_\Omega S^6\right)^\frac{1}{3}\left(\int_\Omega I^3\right)^\frac{2}{3}\leq\left(\int_\Omega S^6\cdot\int_\Omega I^6\cdot|\Omega|\right)^\frac{1}{3}
\]So it follows that $S,I,R\in W^{2,2}(\Omega)$, by the similar argument as above, followed by boostrapping arguments, the claim holds.\vspace{2pt}

\noindent\textbf{Remark:} In the case $n=1$ and $n=3$, we actually proved that the solution is in $C^\infty(\Omega)$, which is a stronger conclusion.

\section{\textcolor{blue}{\textit{Diffusion-driven instability: the Turing instability}}}
The equilibria of the PDE system is important to be used to understand the behavior of the solution.  Define $\bm{A}=(S,I,R)^T$, we can linearize the system around the steady-state $A^*_1$ and $A^*_2$. To do the stability analysis, we linearize the PDE in matrix form\begin{equation}
\partial_t\bm{A}=\bm{J}\bm{A}+\bm{D}\bm{\nabla}^2\bm{A},
\end{equation}where $\bm{J}$ is the differential of the function $(f_1,f_2,f_3)^T$ and $\bm{D}$ is the diagonal matrix consisting of the diffusion constant, i.e. $\bm{D}=\text{diag}(\chi_S,\chi_I,\chi_R)$. Here we may look for a special type of solution to investigate the stability\[
\bm{A}=\eta e^{\lambda t}\bm{B}(\bm{x})\quad \eta\in\mathbb{R}.
\]I.e. the separable solution of $A$ with $\lambda$ a real constant denoting the growth rate around the equilibrium. In addition, we assume $\bm{B}(\bm{x})$ to be the solution to the Homeolz equation\begin{equation}
\bm{\nabla}^2 \bm{B}+k^2\bm{B}=0,\quad \text{and } \partial_n\bm{B}=0 \text{ on $\partial\Omega$},\ k\in\mathbb{R}.
\end{equation}This is motivated by the fact that when the dimension is one, as discussed in [4], the solution will be found to be trigonometric functions, and this will be a Fourier series when writing the sum of separable solution, so $k$ will denote the wave number of the solution. And then, we establish the following result regarding the stability of these equilibria, combining (13) and (14), we see that $\lambda$ satisfies\[
\lambda\bm{B}-\bm{J}\bm{B}-D\bm{\nabla}^2\bm{B}=\left(\lambda\bm{I}-\bm{J}+\bm{D}k^2\right)\bm{B}=0,
\]To guarantee a non-trivial solution, we need $\det(\lambda\bm{I}-\bm{J}+\bm{D}k^2)=0$, i.e.\begin{align*}
&\det\begin{pmatrix}
    \lambda+\beta I+\nu+\chi_Sk^2&\beta S&0\\
    -\beta I&\lambda-\beta S+\gamma+\nu+\chi_I  k^2&0\\
    0&-\gamma&\lambda+\nu+\chi_R k^2
\end{pmatrix}_{\bm{A}^*}=0.
\end{align*}By verification, $\lambda_1=-\nu-\chi_R k^2$ is a root to this cubic equation. Now we consider the equilibria separately. For $\bm{A}^*_1$, the other two roots are found to be\[
\lambda_2=-\chi_S k^2-\frac{\beta b}{\nu}-\nu,\quad \lambda_3=-\chi_I k^2+\frac{\beta b}{\nu}-\gamma-\nu
\]So we see that the critical value (bifurcation) value of the stability of this equilibrium is when $\lambda_3=0$, i.e.\[
k_*^2=\frac{1}{\chi_I}\left(\frac{\beta b}{\nu}-\gamma-\nu\right).
\]Hence the equilibrium is stable when $k^2>k_*^2$ and unstable when $k^2<k_*^2$. Regarding $\bm{A}^*_2$, $\lambda_{2,3}$ can be expressed as\begin{equation*}
\lambda_{2,3}=\frac{-(a+d)\pm\sqrt{(a-d)^2+4bc}}{2},
\end{equation*}where\[
a=\chi_S k^2+\frac{b\beta}{\gamma+\nu},\quad b=\gamma+\nu,\quad c=\beta\left(\frac{b}{\gamma+\nu}-\frac{\nu}{\beta}\right),\quad d=\chi_Ik^2.
\]Obviously the stability depends on the larger eigenvalue, it can be easily seen that the necessary and sufficient condition on the stability is that $ad>bc$. So the critical value satisfies\[
\chi_S\chi_Ik_*^4+\frac{b\beta \chi_I}{\gamma+\nu}k_*^2-\beta b-\nu(\gamma+\nu)=0.
\]
When we are discussing the instability of an equilibrium, there's an interesting phenomenon called Diffusion-driven instability, also known as Turing instability, which was emphasized in [2]. It is when the instability occurs in the presence of the diffusion term "$\bm{\nabla}^2$" and stable in the absence of the diffusion term. Now we consider the following system of ODEs\begin{align*}
    \frac{dS}{dt}&=b-\beta SI-\nu S,\nonumber\\
    \frac{dI}{dt}&=\beta SI-(\gamma+\nu)I,\\
    \frac{dR}{dt}&=\gamma I-\nu R.\nonumber
\end{align*}We aim to study the stability in the case of no diffusion, so we set the $k$ in (14) to be zero, the threshold for stability is also equivalent to set $k=0$ in eigenvalues involving $k$ we just found, thus for the first equilibrium $\bm{A}^*_1$,\[
\lambda_1=-\nu,\quad \lambda_2=-\frac{\beta b}{\nu}-\nu,\quad \lambda_3=\frac{\beta b}{\nu}-\gamma-\nu.
\]The equilibrium is stable when\[
\frac{b}{\gamma+\nu}-\frac{\nu}{\beta}<0.
\]
The first eigenvalue of $\bm{A}^*_2$ is negative, and the other two can be computed by plugging into the value $\beta$. It can be verified easily that the equilibrium is stable when\[
\frac{b}{\gamma+\nu}-\frac{\nu}{\beta}<0.
\]Which is same as the first condition. However, condition this makes the second equilibrium unphysical as we always expect the valid solutions/ equilibria of our PDE model non-negative. So for $\bm{A}^*_2$ to "exists", the equilibrium is unstable. In conclusion, there is Turing instability when $\bm{A}^*_2$ exists and \[
k^2>\frac{-b\beta \chi_I+\sqrt{(b\beta \chi_I)^2+4\chi_S\chi_I(b\beta+\nu\gamma+\nu^2)(\gamma+\nu)^2}}{2\chi_S\chi_I(\gamma+\nu)}.
\]\textbf{Remark:} The Turing instability often leads to Pattern formation in the system. Which can be controlled by choosing suitable domain.
\section{\textcolor{blue}{\textit{Some concluding remarks}}}
In the recent study of the SIR reaction-diffusion model, many useful methods have been developed to analyze the dynamics of the PDE system. For example, in [9], the parameters are replaced by varying functions $\beta(x)$, $\gamma(x)$, etc. So that we can define two different set according to the rate of local transmission of the disease\[
H^-=\{x\in\Omega:\beta(x)<\gamma(x)\}\quad \text{and}\quad H^+=\{x\in\Omega:\beta(x)>\gamma(x)\}.
\]Where the former one denotes the \textit{low-risk cite} and the later one stands for the \textit{high-risk cite}. Unsurprisingly, under certain assumptions on the non-negativity of $S$ and  $I$ with $H^-,H^+\neq \emptyset$, one can investigate the global asymptotic stability of the disease-free equilibrium and endemic equilibrium. This was done by introducing the \textit{basic reproduction number} $\mathcal{R}_0$, which is defined by\[
\mathcal{R}_0=\sup_{\varphi\in H^1(\Omega)\backslash \{0\}}\left\{\frac{\int_\Omega \beta\varphi^2}{\int_\Omega d_I|\bm{\nabla} \varphi|^2+\gamma\varphi^2}\right\}.
\]And the stability results depend on the value of $\mathcal{R}_0$. 
\section{\textcolor{blue}{\textit{References}}}
\noindent [1] Parkinson, C. and Wang, W. (2023). \textit{Analysis of a Reaction-Diffusion SIR Epidemic Model with Noncompliant Behavior}.\par

\noindent [2] Turing, A. (1952). \textit{The chemical basis of morphogenesis}. Philos. Trans. Roy. Soc. B, 237, 37-72.\par

\noindent [3] Ghergu, M. and Radulescu, V. D. (2012). \textit{Nonlinear PDEs: Mathematical Models in Biology, Chemistry and Population Genetics}. Springer.\par

\noindent [4] Peng, R., \& Wang, M. (2005). \textit{Pattern formation in the Brusselator system.} Journal of Mathematical Analysis and Applications, 309(1), 151-166. \par

\noindent [5] Gilbarg, D. and Trudinger, N. S. (2001). \textit{Elliptic Partial Differential Equations of Second Order}. Berlin: Springer.\par

\noindent [6] Shkoller, Steve. \textit{Notes on Lp and Sobolev Spaces.}\par

\noindent [7] Ammar Khanfer. \textit{Applied Functional Analysis}. Springer, 2024.\par

\noindent [8] Moulay Rchid Sidi Ammi, Achraf Zinihi, Aeshah A. Raezah, \& Yassine Sabbar. (2023). \textit{Optimal control of a spatiotemporal SIR model with reaction–diffusion involving p-Laplacian operator.} \textit{Results in Physics, 52}, 106895. \par

\noindent [9] Linda J. S. Allen, B. M. Bolker, Yuan Lou, A. L. Nevai. \textit{Asymptotic profiles of the steady states for an SIS epidemic reaction-diffusion model. Discrete and Continuous Dynamical Systems}, 2008, 21(1): 1-20.\par
\end{document}